\documentclass[12pt,reqno]{article}

\usepackage{amsmath,amssymb,amsfonts,amsthm,mathrsfs,verbatim} 
\usepackage{setspace}
\usepackage{authblk}
\usepackage{caption}
\usepackage{geometry}
\usepackage{graphics} 
\usepackage{graphicx}  
\usepackage{hyperref}
\usepackage{indentfirst}
\usepackage{makeidx}
\usepackage{mathptmx}
\usepackage{sectsty}
\usepackage{color}
\usepackage{bm}

\font\myfont=cmr12 at 18pt

\renewenvironment{abstract}
 {\par\noindent\textbf{\abstractname.}\ \ignorespaces}
 {\par\medskip}

\sectionfont{\centering\fontsize{14}{14}\selectfont\bfseries}
\subsectionfont{\centering\fontsize{13}{13}\selectfont\bfseries}
\subsubsectionfont{\centering\fontsize{12}{12}\normalfont\itshape}

\DeclareCaptionLabelFormat{UC}{\MakeUppercase{#1} #2}
\makeatletter
\renewcommand\@biblabel[1]{#1.}
\makeatother

\newtheorem{theorem}{Theorem}[section]

\theoremstyle{definition}

\newtheorem{algorithm}[theorem]{Algorithm}

\def\E{{\mathbb E}}
\def\R{{\mathbb R}}

\newcommand{\nnb}{\nonumber}

\title{\myfont \bf A Robust Decision Making Framework for Optimal Strategy Selection in Warfare under Model Uncertainty}

\date{}
\author[a,b]{Georgios I. Papayiannis}

\affil[a]{\footnotesize {\it Hellenic Naval Academy, Department of Naval Sciences, Section of Mathematics, Piraeus, GR}}
\affil[b]{\footnotesize {\it Athens University of Economics and Business, Stochastic Modeling and Applications Laboratory, GR }}

\begin{document}
\maketitle

\begin{abstract} In this paper is presented a framework for treating uncertainty in optimal decision problems occuring in combat situations, in order to robustly select the optimal strategy. A stochastic version of the popular Lanchester's aimed-fire model is considered as the underlying combat system describing the combet dynamics, and upon this an optimal decision rule for allocating forces is constructed. This approach results to a very extendable optimal decision framework, where the optimal strategy is chosen by simultaneously treating robustly uncertainty regarding critical combat parameters. 
\end{abstract}


\section{INTRODUCTION}

Optimal decision making in military sciences is traditionally a hot topic and even more in the modern eras, where the decision framework becomes more and more complex due to technological advances and the new characteristics that are introduced in the modern warfare. Classical practices become obsolete and new decision tools needs to be examined and developed to support the officer in charge in better understand the special characteristics and the stochastic nature of the operational problems under study, employing cutting-edge approaches from decision theory and data science. In this perspective, standard approaches from the fields of warfare dynamics modelling and robust decision theory under uncertainty are discussed and combined, constructing an appropriate framework for robust decision making under fuzzy information regarding several aspects of the operational problem under study. 

Warfare dynamics modeling is a very active scientific field and of major importance for military and naval operations. Following the pioneering work of Lanchester \cite{Lan1916} a great number of research papers have been produced investigating more effective modeling approaches for various kinds of warfare in order to represent, more realistically the possible outcomes of a battle situation. Up to now, a lot of effort has been allocated in combat dynamics modelling and numerous interesting approaches have been proposed in the literature either of deterministic or stochastic nature.

In this paper, a framework for robust decision making in combat situations is proposed, taking into account the battle dynamics, where critical combat parameters are subject to uncertainty. For example, consider the case where two sides are involved in a combat, where several critical aspects of the battle are not precisely known by each side, e.g. opponent's attrition rates, initial strengths, forces mixtures, etc, due to incomplete information about opponent's states. Even some of the same's side parameters may not be exactly estimated (like the attrition rate), since their values may depend on battlefield conditions, weather or other non-measurable factors. In such cases, the states of these parameters may be better represented by a probabilistic model. The probability model may be known (determined by the decision maker/officer in charge), however in many cases is provided by an expert/consultant or a multitude of models are provided by various experts/consultants of varying validity. Then, the decision maker has to allocate her/his preferences on the provided models, i.e. the level of aversion from the provided opinion(s).

In this direction, several aspects of the decision problem has to be determined. First of all, the decision criterion, or the so called performance measure, for the problem has to be determined and of course the control variables of the problem, as long as a certain evolution law (either stochastic or deterministic) for the combat has been adopted. Secondly, if various models about the quantity under question are available, an appropriate aggregation scheme to fuse the available information has to be employed resulting to a single aggregate model which should robustly represent the received information about the unknown parameters of the problem. At last step, all the above are combined to define an optimal control problem concerning the strategy selection, under worst-case scenarios, and numerical schemes for approximating tha optimal solution are proposed.

\section{LITERATURE OVERVIEW \& KEY CONCEPTS}

In this section a brief literature overview on the latest advances on warfare modeling is presented while some key concepts on robust decision making and risk quantification are also discussed. 

\subsection{A Brief Literature Overview On Warfare Dynamics Modelling}%

Combat modelling is a field which has received a lot of attention following the seminal work of Lanchester \cite{Lan1916}. This classic (however still valid in many occasions) model, known as {\it Lanchester's aimed-fire combat model}, describes the battle conflict between two military forces, let us say the blue ones (B) and the red ones (R). Assume that $B(t), R(t)$ denote the strength at a certaint time instant $t$ of the blue and red side respectively while $r,b>0$ are some constant parameters representing the attrition rate of each side to the opponent. In this case, the battle dynamics are represented by the system of differential equations
\begin{eqnarray}\label{lan-mod}
\left\{ \begin{array}{c}
dB(t) = - r R(t) dt\\
dR(t) = - b B(t) dt
\end{array} \right.
\end{eqnarray}
Beyond its simplicity, this modeling approach has the ability to represent on a very straightforward manner the evolution of a battle. Of course the above approach can be very easily modified for mixed forces battle occations.\\ 

Probably the most famous extension is the Bracken's generalized model \cite{Bracken1995} where some exponential parameters are used in order to capture nonlinear characteristics of the combat, however making the modeling approach more data-dependent. In this spirit, Bracken proposed the model
\begin{eqnarray}\label{bra-mod}
\left\{ \begin{array}{c}
dB(t) = - r R(t)^p B(t)^q dt\\
dR(t) = - b B(t)^p R(t)^q dt
\end{array} \right.
\end{eqnarray}
where the exponential parameters $p,q>0$ need to be determined (either emprically or using experts' opinions) allowing the contribution of a scaled by $p,q$ interaction of both forces to each side's attrition rate. The classic Lanchester model can be retrieved for $q=0$ and $p=1$.\\ 

\begin{figure}[ht!]
\centering
    \begin{minipage}{.48\textwidth}
        \centering
        \includegraphics[width=3in]{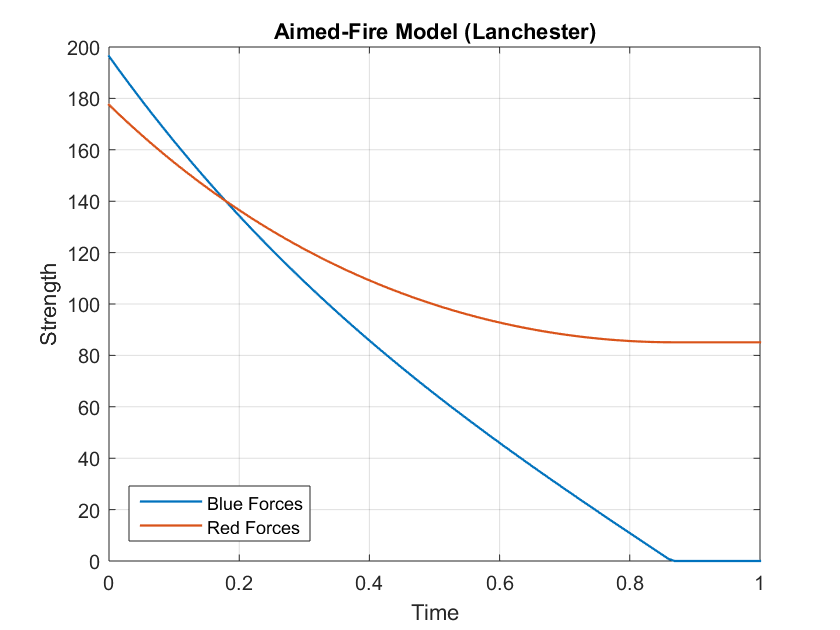}
    \end{minipage}%
    \begin{minipage}{.48\textwidth}
        \centering
        \includegraphics[width=3in]{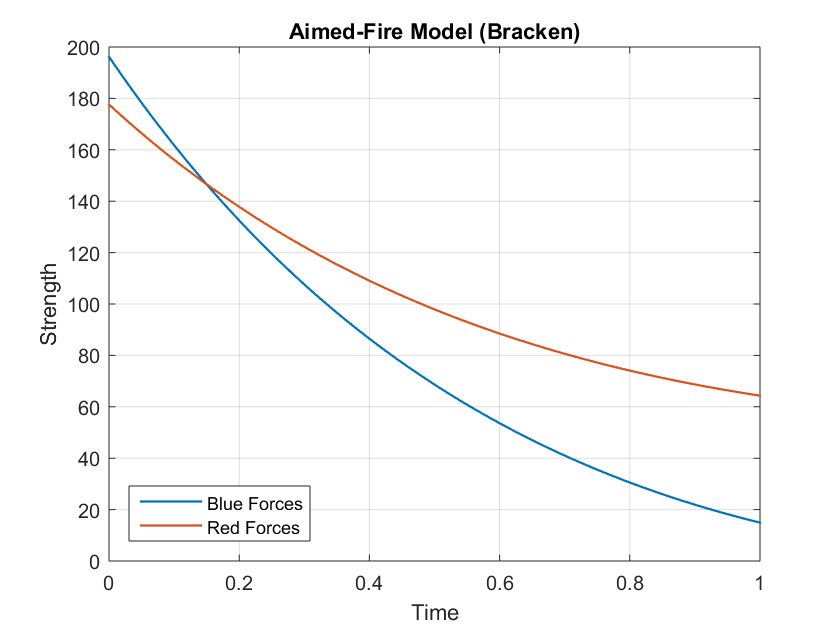}
    \end{minipage}%
    \caption{Warfare trajectories generated by Lanchester's model (left column) and Bracken's model (right column).}\label{fig-1}
\end{figure}

Some very interesting results from a stochastic perspective were shown in \cite{King2002} where the classic Lanchester's model is reformulated to 
\begin{eqnarray}\label{king-mod}
\left\{ \begin{array}{c}
dB(t) = - r( R(t) )\,\, h(B(t), R(t), t) dt\\
dR(t) = - b( R(t) )\,\, h(B(t), R(t), t) dt
\end{array} \right.
\end{eqnarray}
with $r(\cdot), b(\cdot), h(\cdot, \cdot, t)$ being all strictly positive functions. This approach allows the firepower of each army (or fleet) to depend nonlinearly on the number of the forces (or the ships) where battlefield and weather conditions can be captured through an appropriate choice of the function $h$. This model has been examined under a discrete Markovian framework, however this formulation offers a great flexibility to the represention of different warfare dynamics patterns. More recently, \cite{kim2017a} studied the quite simple but very interesting stochastic version of the aimed-fire model
\begin{eqnarray}\label{kim-mod-a}
\left\{ \begin{array}{cc}
dB(t) = - X r R(t) dt, & X \sim Bernoulli(p) \\
dR(t) = - Y b B(t) dt, & Y \sim Bernoulli(q) 
\end{array} \right.
\end{eqnarray}
where $X,Y$ Bernoulli random variables with probability $p$ that attrition rate is realized (either for both sides or to one of them). In Figures \ref{fig-1} and \ref{fig-2} warfare dynamics trajectories for the aimed-fire model under the deterministic and the stochastic framework are illustrated.

\begin{figure}[ht!]
\centering
    \begin{minipage}{.48\textwidth}
        \centering
        \includegraphics[width=3in]{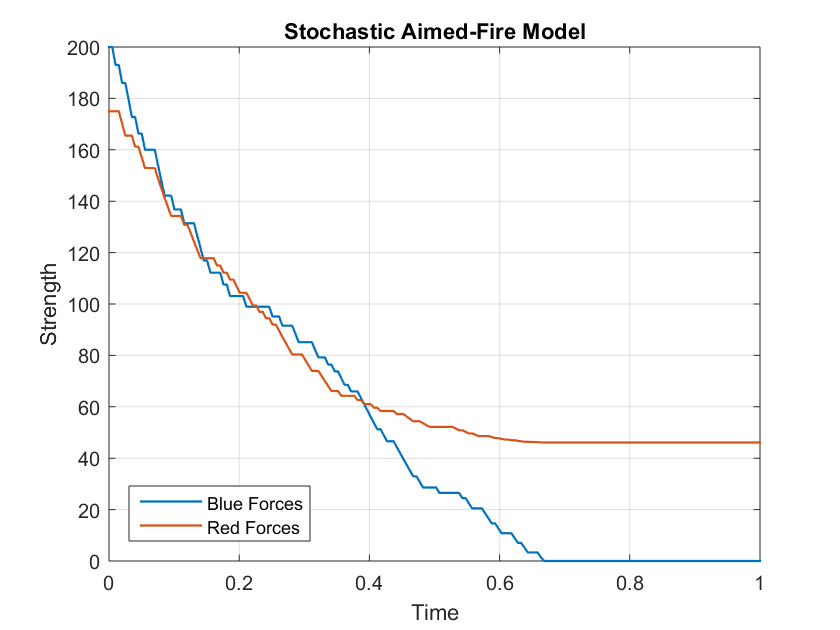}
    \end{minipage}%
    \begin{minipage}{.48\textwidth}
        \centering
        \includegraphics[width=3in]{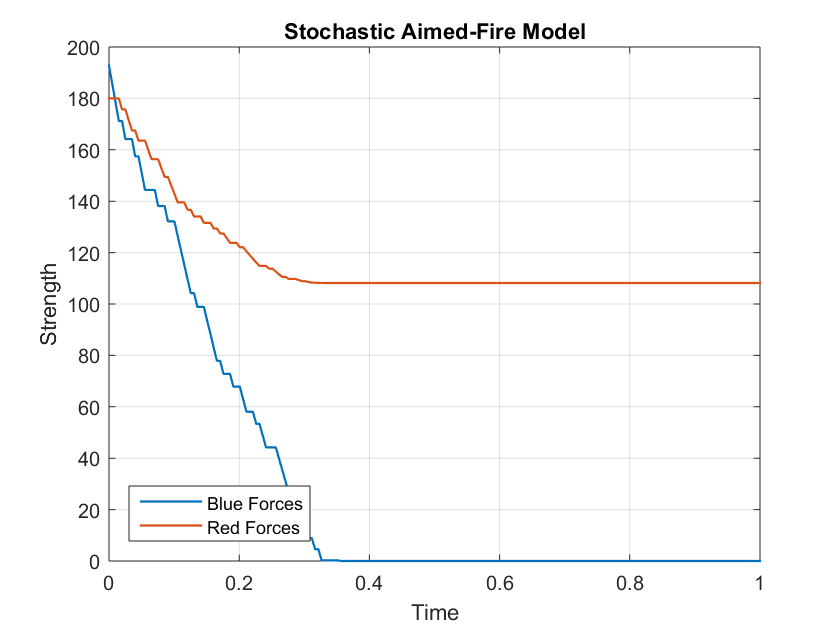}
    \end{minipage}%
    \caption{Warfare trajectories simulated by stochastic aimed-fire model using independence assumption between the success probabilities in attrition rates realization (left column) and strong negative dependence (right column).}
\end{figure}\label{fig-2}

Finally in \cite{kim2017b} the aimed-fire model is examined under the framework of the network centric warfare (NCW) where one of the sides (let us say the red ones) is assumed to be supported by a network $N$. In this case the combat dynamics (for a network with two agents supporting the red side) are described by the system
\begin{eqnarray}\label{kim-mod-b}
\left\{ \begin{array}{l}
dB(t) = - g(N) R(t) dt\\
dR(t) = - \pi_0 b_R B(t) dt\\
dA^1(t) = -\pi_1 b_1 B(t) dt\\
dA^2(t) = -\pi_2 b_2 B(t) dt
\end{array} \right.
\end{eqnarray}
where $b_R,b_1,b_2$ denote the attrition rates of the blue side (B) to the red side directly (R) and their network's agents $A^1, A^2$ respectively; $\pi_0, \pi_1, \pi_2$ the blue side's fire allocation proportion directly to the red side and the supporting agents and $g(N)$ the attrition function of the red side, where the attrition rate depends on the network connectivity. The latter becomes much more complicated when networks on both sides are considered however is one of the most realistic models for modern warfare. Although we do not examine this case in this work, the proposed decision making framework is directly applicable to such systems of warfare dynamics uner minor modifications.


\subsection{Model Uncertainty, Opinions Aggregation and Robust Information Treatment}\label{sec-2.2}%

{\it Model uncertainty} or {\it model ambiguity} is a key issue in modern decision making. In general, this term refers to the ignorance of the decision maker regarding a plausible model that efficiently describes certain parameters of the problem under study. The typical setting is that decision maker exchanges her/his ignorance to some expert opinions on the matter that she/he cannot provide a model by herself/himself. As long as a model is provided, the decision maker has of course the choice to deviate from it, if the model does not seem realistic or if no evidence about the reliability of the expert is available, leading to a worst-case decision framework \cite{Gilboa1989}. Assume that the decision maker's preferences are quantified by a utility function $U(\alpha, Z): \mathcal{A} \times \mathcal{Z} \to \R$, where $\mathcal{A}$ denotes the action set (and $\alpha$ the action), $\mathcal{Z}$ denotes the set of the unknown problem parameters (with $Z$ denoting a particular scenario) and $Q_0$ the model provided by the expert for the possible states of $Z$. Then, robust decision rules for optimal selection of $\alpha \in \mathcal{A}$ can be derived through the {\it multiplier preferences} approach \cite{Mac2006}
\begin{equation}\label{multiplier-preferences}
V(\alpha) = \min_{ Q \in \mathcal{P}(\mathcal{Z}) } \left\{ \mathbb{E}_Q[U(\alpha, Z)] + \theta d(Q,Q_0) \right\}, \,\,\, \theta>0
\end{equation}
where parameter $\theta>0$ quantifies the decision maker's aversion preferences from $Q_0$ and $d$ is an appropriate metric in the space of probability models. Equivalently, the same decision rule can be represented through the {\it constraint preferences} approach \cite{Hansen2011}  
\begin{equation}\label{constraint-preferences}
\widetilde{V}(\alpha) = \min_{ Q \in \mathfrak{Q}_c } \mathbb{E}_Q[U(\alpha, Z)] 
\end{equation}
where $\mathfrak{Q}_c := \{ Q\in\mathcal{P}(\mathcal{Z})\,\, : \,\, d(Q,Q_0)\leq c \}$ denotes the set of plausible models for $Z$ with respect to the decision maker's aversion preferences from $Q_0$ which are quantified by the parameter $c>0$. It is quite simple to show that for certain choices of parameters $\theta>0$ and $c>0$, given that the same metric sense $d$ is used, both approaches coincide. This means that if for a certain choice $\theta_0>0$ problem \eqref{multiplier-preferences} is solved leading to an evaluation measure $Q_* := Q_{\theta_0}$, then there exists a certain choice $c_0>0$ for which the problem \eqref{constraint-preferences} leads to an evaluation measure $Q^* := Q_{c_0}$ holding that $Q_* \equiv Q^*$ resulting to the same valuation rule for every possible action $\alpha \in \mathcal{A}$.

Clearly, an important matter on the aforementioned rules is the choice of the metric $d$ which is used to quantify divergence/dissimilarity between the various models. A typical metric sense that is used in the decision theory framework (but also in the information sciences in general) is the Kullback-Leibler divergence \cite{KL1951}, or else known as relative entropy, defined as 
\begin{equation}\label{KLD}
KL(Q \| Q_0) := \int_{\mathcal{Z}} f(z) \log \frac{ f(z) }{ f_0(z) } dz
\end{equation} 
where $f, f_0$ denote the respective probability density functions of the probability models $Q,Q_0$. A very important matter on model uncertainty treatment arises when more than one expert opinions (probability models) are provided to the decision maker. In this case, the opinion aggregation is of particular importance in order to properly take into account all different aspects of information that are available in the most efficient manner. This task depends strongly on the metric choice $d$ that is chosen to provide an appropriate model $Q$ for the description of the stochastic factors of the problem in question. Assume that $N$ different expert opinions are provided through the set $\mathcal{Q} := \{Q_1, Q_2,..., Q_n\}$. Generalizing the approaches presented in \eqref{multiplier-preferences} and \eqref{constraint-preferences}, it is needed to express the distance of a particular probability model $Q$ from a whole set of probability models like $\mathcal{Q}$. Under this {\it multi-prior models} framework, modified versions of metrics in the space of probability measures are employed. In this fashion, the multi-prior extension of the Kullback-Leibler divergence (or KL-barycenter) is defined as
\begin{equation}\label{w-KLD}
KL_w(Q \| \mathcal{Q}) := \sum_{i=1}^N w_i \int_{\mathcal{Z}} f(z) \log \frac{ f(z) }{ f_i(z) } dz
\end{equation} 
where $f_i$ denote the respective probability density function of the probability model $Q_i$ for $i=1,2,...,N$ and the weights $w$ are such that $\sum_{i=1}^N w_i=1$ and $w_i \geq 0$ for all $i$. The weights could be either allocated by the decision maker representing her/his beliefs regarding each prior model's validity or, if not such information is available, could be all allocated uniformly $w_1 = ... = w_N = \frac{1}{N}$. Obtaining the minimizers of such weighted versions of metrics is equivalent to averaging the provided by the experts models with respect to the sense of the metric is used, e.g. minimizing \eqref{w-KLD} with respect to $Q$, corresponds to estimating the Kullback-Leibler barycenter of the set $\mathcal{Q}$ leading to the geometric mean of the set. In the case where no aversion preferences from the prior-models set are allocated by the decision maker, then the barycenter is the best possible representation of the information contained in the set $\mathcal{Q}$ and under this model the decision rules in \eqref{multiplier-preferences} and \eqref{constraint-preferences} are evaluated. If certain aversion-preferences are determined, then distorted versions of the respective barycenters are the optimal representations of the priors' set and under these the optimal actions of the decision maker are chosen. For more details on belief aggregation issues and some applications in policy selection, actuarial and financial applications please see \cite{Papayiannis2018} and \cite{Petracou2022}.

\section{ROBUST DECISION MAKING UNDER STOCHASTIC WARFARE DYNAMICS DRIVEN BY MODEL UNCERTAINTY}

In this section, a framework for robust decision making under uncertainty in combat situations is proposed when stochastic warfare dynamics are considered. As an illustrative example is discussed the stochastic aimed-fire model.

\subsection{A General Robust Decision Framework}\label{sec-3.1}

Assume that two sides participate in a conflict, the blue ones (B) and the red ones (R). Let $\pi$ denote the control (strategy, e.g. proportion of available weapons/forces to be used) of the blue side, i.e. from now on the states of $B(t), R(t)$ depend both on $\pi$ so the notation is changed to $B^{\pi}(t), R^{\pi}(t)$ in order to clarify this dependence. Assume further that the dynamics of the battle can be efficiently calibrated by the stochastic version of the aimed-fire model resulting to the system
\begin{eqnarray}\label{mod-0}
\left\{ \begin{array}{c}
dB^{\pi}(t) = - r(Z) R^{\pi}(t) dt,\,\,\, B^{\pi}(0)=B_0, \,\,\, Z_B(t) : \Omega \to \R^{n_B} \\
dR^{\pi}(t) = - b(Z) B^{\pi}(t) dt,\,\,\, R^{\pi}(0)=R_0, \,\,\, Z_R(t) : \Omega \to \R^{n_R}
\end{array} \right. , \,\,\, Z = (Z_B , Z_R) \sim Q, 
\end{eqnarray}
where $B_0,R_0$ are the initial states (initial strengths) for both sides, $\Omega$ denotes the set of all possible situations that may occur on the battlefield that could affect the attrition rate of each side, e.g. weather or terrain conditions, effects of external forces, etc., $Z_B, Z_R$ the stochastic factors for both sides related to events in $\Omega$ that directly affect each side's attrition rate, $Q$ denotes the law that describes the random behaviour of both random variables $Z_B, Z_R$ and $b:\R^{n_B}\times R^{n_R}\to \R$, $r:\R^{n_B}\times R^{n_R}\to \R$ the attrition rate mappings for each side. Note here that no knowledge on opponent's (red side) strategy is assumed and under this perspective a worst-case scenario is adopted under which all red forces are assumed to be employed into the battlefield to achieve their maximum strength and attrition rate to the blue side. This is the case of worst-case framework where everything that is uknown is considered as the force of the nature and the nature is assumed to play totally against the decision's maker side. 

First, the decision maker is interested to the probability model $Q \in \mathcal{P}(\Omega)$ that could best describe the evolution of the battle in order to simulate in the most realistic manner the scenarios that may occur. Under the model uncertainty concept, the decision maker (blue side) is provided with a set of $N$ plausible probability models, let us denote it by $\mathcal{Q}$, from various battle experts or information sources. For each one of them, a certain trustworthyness weight $w_i$ is allocated by the decision maker. A strategy $\pi$ needs to be choden in order to optimize certain criteria and possibly on a certain time interval, let us denote it by $[0, T].$ The criteria taken into account could be: the expected state of the blue forces with respect to the red ones at a terminal time instant $T$: $\E_Q[B(T; \pi)-R(T)]$, the expected total casualties of the blue forces $\E_Q[B(T;\pi)-B(0;\pi)]$, the value of the reserved forces at terminal time $\E_Q[\psi(B(T;\pi), R(T) )]$ where $\psi$ the function evaluating the value of reserves at $T$ depending on the final state of $B(t),R(t)$, etc. However, the aforementioned criteria depend on the probability model $Q$ highlighting the importance and how crucial is the choice of this ``valuation measure'' in the selection of the optimal decision. 

Following the discussion in Section \ref{sec-2.2}, for any appropriate utility function $\varphi$ depending on $\pi$ and $Z$, the dynamics of $B(t)$ and $R(t)$ as described in \eqref{mod-0}, given decision maker's aversion preferences from the prior models set $\mathcal{Q}$ (either in terms of \eqref{multiplier-preferences} or \eqref{constraint-preferences}) and credibility allocation to each model in the set, a generic formulation of the strategy selection problem is achieved as follows:
\begin{eqnarray}\label{robust-DM}
&&\max_{\pi \in \mathcal{A}} \min_{Q \in \mathfrak{Q} } \Phi(\pi, Q) = \max_{\pi \in \mathcal{A}} \left\{ \min_{Q \in \mathfrak{Q} }\E_Q[\varphi(Z(T); \pi)] \right\} \\
&&\mbox{subject to}\nnb \\
&&(B(t), R(t)) \mbox{ described by dynamical system in \eqref{mod-0}} \nnb
\end{eqnarray}
where $\mathcal{A}$ contains all admissible strategies $\pi$, $\mathfrak{Q}$ is the set of available probability models describing the dynamics of $B,R$ satisfying the aversion preferences of the decision maker from the available information (prior set $\mathcal{Q}$). 

In the case where Kullback-Leibler divergence is used as metric tool as already discussed, the inner's problem minimizer can be explicitly calculated for any $\pi$, therefore the above problem reduces to a maximization problem with respect to $\pi$ since the minimizer $Q_*$ can be obtained in closed form. In particular, we refer to the result from \cite{Papayiannis2018} where given a prior set $\mathcal{Q}$ of $N$ provided models, validity weights $w\in\Delta^{N-1}$ and aversion preferences $\theta>0$ (i.e. in terms of \eqref{multiplier-preferences}), the minimizer $Q_\theta$ of the inner problem stated in \eqref{robust-DM} is explicitly represented (in terms of the respective probability density function) by
\begin{equation}\label{Qtheta}
f_{\theta}(z) = C_{\theta}\,\, e^{ -\theta \varphi(z) }\,\, \prod_{i=1}^N f_i(z)^{w_i}, 
\end{equation}
where $C_{\theta}^{-1} := \int e^{ -\theta \varphi(z) } \prod_{i=1}^N f_i(z)^{w_i}dz$. Importantly, in the limit case where $\theta \to \infty$ (absolute trust in the prior set is allocated by the decision maker) the minimizer $Q_*$ coincides with the KL-barycenter, i.e. 
\begin{equation}\label{Qbar}
f_{*}(z) = C_{*}\,\,\prod_{i=1}^N f_i(z)^{w_i}, 
\end{equation}
where $C_{*}^{-1} := \int \prod_{i=1}^N f_i(z)^{w_i}dz$. Therefore plugging \eqref{Qtheta} to \eqref{robust-DM} the optimization problem becomes very computationally tractable since from a minimax problem is reduced to a maximization problem with respect to $\pi$. In the next subsection, a discretized version of the Bernoulli-distributed stochastic aimed-fire model is discussed in detail, simulation schemes for the dynamics are proposed and a stochastic gradient approach for the optimal solution of the resulting control problem is described. 

\subsection{Optimal Strategy Selection under a Stochastic Aimed-Fire Model}\label{sec-3.2}%

As a concrete example, an optimal strategy selection problem is considered under a stochastic version of the Lanchester model presented in \eqref{kim-mod-a} under the framework of model uncertainty. 

\subsubsection{A Discretized Stochastic Version of the Aimed-Fire Model}%

Consider that the blue's side commander (decision maker) desires to statically select up to a certain time instant $T$ (where the battle conditions are expected to change) the most efficient allocation for the available forces. Let us denote by $\pi \in (0,1]$ the proportion of blue forces to be used in the battle and by $1-\pi$ the proportion of forces to be kept out of the battle in order to be used after time $T$. Red side is assumed to use all of its available power since a worst-case scenario strategy is going to be derived here for the blue side. 

We consider a discretized version of the stochastic combat model \eqref{kim-mod-a} and in particular we consider the model
\begin{eqnarray}\label{sd-1}
\left\{ \begin{array}{ll}
B^{\pi}_{k \Delta t} = B^{\pi}_{(k-1)\Delta t} - r \, Z_{k,R} R^{\pi}_{(k-1) \Delta t}, & B^{\pi}_0 = B_0, \,\,\, k=1,2,... \\
R^{\pi}_{k \Delta t} = R^{\pi}_{(k-1)\Delta t} - \pi\, b\, Z_{k,B} B^{\pi}_{(k-1) \Delta t}, & R^{\pi}_0=R_0, \,\,\, k=1,2,...\\
Z_k = (Z_{k,B}, Z_{k,R}) \sim Q &
\end{array} \right. 
\end{eqnarray}
where $Z_k$ for $k=1,2,...$ are considered as independent and identically distributed random variables according to a probability law $Q$ and $\Delta t >0$ is considered as a relative small time increment. In order to provide a more specific example, assume that $Q$ is a probability model with Bernoulli marginals, i.e.
$$ Z_{k,B} \sim Bernoulli(p_B), \,\,\,\, Z_{k,R} \sim Bernoulli(p_R) $$
with $p_B, p_R$ denoting the probabilities of attrition realization for each side, binded by a specific dependence relation, represented by some copula $C$. Note that the model of \cite{kim2017a} considers only a particular case of the above model, i.e. the case where the independence copula is used. If other than the independence copula is used, i.e. assuming some kind of dependence between the probabilities of attrition realization, can offer a more general framework where various combat conditions and dimensions may be captured. In this case, if $F$ denotes the joint distribution function of model $Q$ and $F_B, F_R$ denotes the Bernoulli marginals distribution functions, then under celebrated Sklar's theorem \cite{Sklar1973} the distribution function of $Z$ can be represented through the relation
\begin{equation}
F(Z_B, Z_R) = C( F_B(Z_B), F_R(Z_R) )
\end{equation}
for a certain copula function $C:[0,1]\times[0,1] \to [0,1]$ (see \cite{Nelsen2007} for more details on copulae). 

Since the system of state equation presented in \eqref{sd-1} is of homogeneous form, we may derive an explicit solution for each time increment of length $\Delta t$ under the assumption that the realization of the random variable $Z$ is known. Writing the system in \eqref{sd-1} in a matrix formulation we get the linear system
\begin{eqnarray}\label{sol}
X^{\pi}_{k \Delta t} := 
\begin{pmatrix}
B^{\pi}_{k \Delta t}\\ R^{\pi}_{k \Delta t}
\end{pmatrix}
= 
\begin{bmatrix}
0 						  & - r\, Z_{k,R} \, \Delta t\\
-\pi \, b \, Z_{k,B} \, \Delta t & 0 
\end{bmatrix} \cdot 
\begin{pmatrix}
B^{\pi}_{(k-1) \Delta t}\\ R^{\pi}_{(k-1) \Delta t}
\end{pmatrix}
=:
A^{\pi}_k \cdot X^{\pi}_{(k-1) \Delta t}, \,\,\, k=1,2,...
\end{eqnarray}
where $X^{\pi}$ the state variable and $A^{\pi}$ the coefficient matrix of the state system which depends directly on the realization of the random variable $Z$. However, under the probability model $Q$ considered in this variation of the aimed-fire model, only four possible realizations for $Z$ are possible at each step $k$, therefore the matrix $A$ can have only four outputs depending on the variable's $Z = (Z_B, Z_R)$ outcome, i.e. 
\begin{equation}
A_k^{\pi} = \left\{
\begin{array}{ll}
\begin{bmatrix}
0 & 0\\
0 & 0 
\end{bmatrix}, & Z_{k,B} =0, \,\,\, Z_{k,R} = 0\\
&\\
\begin{bmatrix}
0 & 0\\
-\pi b \Delta t & 0 
\end{bmatrix}, & Z_{k,B} =1, \,\,\, Z_{k,R} = 0\\ 
&\\
\begin{bmatrix}
0 & -r \Delta t\\
0 & 0 
\end{bmatrix}, & Z_{k,B} =0, \,\,\, Z_{k,R} = 1\\
&\\
\begin{bmatrix}
0 & -r \Delta t\\
-\pi b \Delta t & 0 
\end{bmatrix}, & Z_{k,B} =1, \,\,\, Z_{k,R} = 1\\
\end{array}
\right.
\end{equation}
for each $k=1,2,...$. Assuming that at each step the realizations of the random variable $Z$ are known (i.e. the occuring matrices $A_k$ are known) then the state equation can be represented in terms of the initial state $X_0$ by applying the analytic solution formula \eqref{sol} recursively for each time step $\Delta t$ and get:
\begin{equation}\label{cf-sol}
\boxed{X^{\pi}_{k\Delta t} = 
    e^{ \Delta t \, A_{k}^{\pi} } \cdot X^{\pi}_{(k-1)\Delta t} 
= e^{ \Delta t \, \sum_{i=1}^k A_{i}^{\pi} } \cdot X_0, \,\,\, k=1,2,...}
\end{equation}

\subsubsection{A Gaussian Approximation Scheme for the Distribution of Combat Dynamics and a Simulation Approach}

For simulation purposes, a Gaussian approximation scheme is described here for the states of the system $X=(B,R)$ at any time step $k$ conditional on the states of $X$ at previous step $k-1$. Under this approximation, it is assumed that 
\begin{equation}
X^{\pi}_{k \Delta t}|X^{\pi}_{(k-1)\Delta t} \sim N( \mu_k, S_k ), \,\,\, k=1,2,..., \,\,\, \mbox{for} \,\,\, \Delta t>0,
\end{equation}
where $\mu_k \in \R^2$ is the vector with the location and $S_k \in \R^{2\times 2}_{+}$ is the positive-definite matrix containing the dispersion parameters. Below up to second-order moments are calculated in order to obtain estimates for the Gaussian distribution parameters.

\begin{eqnarray}
m_{k,B} &=& \E[B^{\pi}_{k\Delta t}|X^{\pi}_{(k-1)\Delta t}] = B^{\pi}_{(k-1)\Delta t} - r R^{\pi}_{(k-1)\Delta t} \E[Z_{k,R}] \Delta_t  \\
m_{k,R} &=&\E[R^{\pi}_{k \Delta t}|X^{\pi}_{(k-1)\Delta t}] = R^{\pi}_{(k-1)\Delta t} - \pi b B^{\pi}_{(k-1)\Delta t} \E[Z_{k,B}] \Delta_t \\
m^2_{k,B} &=& \E[B^{\pi,2}_{t+\Delta t}|X_{(k-1)\Delta t}] = B^{\pi,2}_{(k-1)\Delta t} - 2 r B^{\pi}_{(k-1)\Delta t} R^{\pi}_{(k-1)\Delta t}  \E[Z_{k,B}] \Delta_t \nnb \\
&\,\,\,& + r^2 R_{(k-1)\Delta t}^{\pi,2} \E[Z_{k,R}^2] \Delta^2_t\\
m^2_{k,R} &=&\E[R^{\pi,2}_{k \Delta t}|X^{\pi}_{(k-1)\Delta t}] = R^{\pi,2}_{(k-1)\Delta t} - 2 \pi b B^{\pi}_{(k-1)\Delta t} R^{\pi}_{(k-1)\Delta t}  \E[Z_{k,R}] \Delta_t \nnb \\
&\,\,\,& + \pi^2 b^2 B^{\pi,2}_{(k-1)\Delta t} \E[Z_{k,B}^2] \Delta^2_t\\
m_{k, BR} &=& \E[B^{\pi}_{k\Delta t} R^{\pi}_{k\Delta t}|X^{\pi}_{(k-1)\Delta t}] = B^{\pi}_{(k-1)\Delta t} R^{\pi}_{(k-1)\Delta t} - \pi b B^{\pi,2}_{(k-1)\Delta t}  \E[Z_{k,B}] \Delta t - r R^{\pi,2}_{(k-1)\Delta t}  \E[Z_{k,R}] \Delta t \nnb\\
&\,\,\,& + \pi b r B^{\pi}_{(k-1)\Delta t} R^{\pi}_{(k-1)\Delta t}  \E[ Z_{k,B} Z_{k,R} ] \Delta^2_t
\end{eqnarray}
On the above calculations under the Bernoulli marginals hypothesis we get
$$ \E[Z_{k,B}] = \E[Z^2_{k,B}] = p_B, \,\,\,\,\,\, \E[Z_{k,R}] = \E[Z^2_{k,R}] = p_R$$
while for the covariance term one has to calculate
$$ \E[Z_{k,B}Z_{k,R}] = \int z_{B} z_{R} dC(F_B(z_B), F_R(z_R)) $$
which clearly depends on the copula choice. For instance, if the independece copula is selected $C(u,v) = u\cdot v$, then the covariance term will be $ \E[Z_{k,B}Z_{k,R}] = \E[Z_{k,B}]\cdot \E[Z_{k,B}]=p_B \cdot p_R$. Using the moment calculations, the Gaussian distribution parameters are estimated as 
\begin{eqnarray}\label{gauss-est}
\mu_k = \begin{pmatrix} m_{k,B}\\ m_{k,R} \end{pmatrix}, \,\,\,\,\,\,
S_k = \begin{bmatrix}
m_{k,B}^2 - (m_{k,B})^2 & m_{k,BR}-m_{k,B} m_{k,R}\\
m_{k,BR}-m_{k,B} m_{k,R} & m_{k,R}^2 - (m_{k,R})^2
\end{bmatrix}, \,\,\,\,\,\, k=1,2,...
\end{eqnarray}

Below an algorithm for simulating the states of combat dynamics described in \eqref{sd-1} for given $\pi \in (0,1]$ under the Gaussian approximation scheme is described in steps.

\begin{algorithm}[{\it Aimed-Fire Model Paths Simulationn Scheme}]\label{alg-0}
For a given strategy $\pi \in (0,1]$, time step $\Delta t$ and time horizon $T$ ($N=T/\Delta t$), initial states $X_0 = (B_0, R_0)$, and a dependence structure $C$ (copula), set $k=1$ and repeat the following steps for $k=1,2,...,N$: 
\begin{enumerate}
\item[{\bf 1.}] According to the (known) states $X_{(k-1)\Delta t}$ calculate the Gaussian distribution parameters in \eqref{gauss-est} and denote the estimated probability model by $Q_k$.
\item[{\bf 2.}] Simulate the random states at step $k$ by drawing values $X_{k\Delta t}=(B_{k\Delta t}, R_{k \Delta t}) \sim Q_k$.
\item[{\bf 3.}] Set $k = k+1$.
\end{enumerate}
\end{algorithm}

\subsubsection{The Optimality Criterion and the Optimal Strategy Selection Problem}%

Assume that the commander of the blue side (decision maker) desires to take into account a number of criteria while choosing his best strategy $\pi$. Clearly, optimality of the selected strategy depends directly on the preferences of the decision maker regarding the level of importance (weights) that allocates in the various criteria that are taken into account. For the sake of brevity, assume that three criteria should be taken into account for the selection of strategy $\pi$: (a) difference from the oppositions strength at the final time instant $T = N \Delta t$, (b) difference from a specified minimum level of strength $B_{\min} \leq B_0$ at the final time instant $T$ and (c) available fresh forces to be used after time instant $T$ where the battle conditions are expected to change.

A standard technique in decision theory is to aggregate all required criteria to a single objective function by introducing some sensitivity parameters assigned to each criterion representing the importance allocated by the decision maker to each one of the criteria. In particular, a scaled version is used here where each parameter will indicate the proportion of importance allocated to each particular criterion taking also into account possible scale effects in the objective function. Therefore, the payoff/utility function that is considered for the problem under study is
\begin{equation}\label{utility-1}
\Phi(\pi) := \theta_1 \E_Q [ \varphi(B_T^{\pi}-R_T^{\pi}) ] + \theta_2 \E_Q[ \varphi(B_T^{\pi} - B_{\min}) ] + \theta_3 \psi( (1-\pi) B_0), \,\,\, \sum_{i=1}^3 \theta_i=1, \,\, \theta_i \geq 0 
\end{equation}
where $Q$ denotes the probability model under which the states $B,R$ evolve, $\varphi(u) := u^2 (I_{u\geq 0}(u) - I_{u<0}(u))$ is a quadratic profit function used to measure deviance of the blue side's strength from certain targets, $\psi(u) := e^{\zeta T}u$ denotes the importance function for keeping an ammount of strength out of the combat up to the terminal time $T$ weighted by a sensitivity parameter $\zeta>0$ which is suitably chosen either from the decision maker (according to her/his preferences) or the combat conditions and what is expected to happen after time $T$ on the battlefield. The parameter vector $\theta = (\theta_1, \theta_2, \theta_3)$ allocate the decision maker's importance to three criteria described above, and specifically: (a) $\theta_1>0$ quantifies the importance allocated in achieving the greatest positive difference from the opponents strength u pto time instant $T$, (b) $\theta_2>0$ represents the importance in keeping available forces strength higher than a certain thresshold $B_{\min}<B_0$ up to time instant $T$ and (c) $\theta_3>0$ is the importance placed in keeping a proportion $1-\pi$ of the available power out of the combat up to terminal time $T$. In this manner, various criteria are aggregated in a single payoff function, however the final decision criterion relies strongly on the choice of the sensitivity and weighting parameters $\theta, \zeta$ which are a priori defined.

In this fashion, the optimal strategy selection problem can be stated as
\begin{eqnarray}\label{opt-0}
&&\max_{\pi \in [0,1]} \Phi(\pi) \\
&&\mbox{subject to}\nnb \\
&& \left\{ \begin{array}{ll}
B^{\pi}_{k \Delta t} = B^{\pi}_{(k-1) \Delta t} - r Z_{k,R} R^{\pi}_{(k-1)\Delta t} \Delta t, \,\,\, k=1,2,...,N, & B^{\pi}_0 = B_0 \\
R^{\pi}_{k \Delta t} = R^{\pi}_{(k-1) \Delta t} - \pi b Z_{k,B} B^{\pi}_{(k-1)\Delta t} \Delta t, \,\,\, k=1,2,...,N,& R^{\pi}_0=R_0\\
Z_k = (Z_{k,B}, Z_{k,R}) \sim Q,\,\,\, \forall k=1,2,...,N &
\end{array} \right. \nnb
\end{eqnarray}
where probability model $Q$ is considered known.

\subsubsection{The Model Uncertanty Issue and the Proposed Framework}

In most cases, exact knowledge of model $Q$ is impossible or a multitude of models are provided as different experts' opinions. This fact causes the problem of model uncertainty which has to be simultaneously treated with the problem of optimal strategy selection. In the case considered above, assume that there is model uncertainty regarding parts (opponent's stochastic parameters and/or correlation structure) or the whole probability model $Q$ (consider also fuzzy the parameter $p_B$ due to incomplete knowledge of the opponent's defense abilities). In this case, the problem's \eqref{opt-0} objective function should be substituted by $J(\pi, Q)$ depending both on the strategy $\pi$ and the probability model $Q$. A property that maybe required is that $ \Phi(\pi) = J(\pi, Q^*)$ where $Q^*$ is the optimal probability model, i.e. the true one. However, the approximation of the true probability law is a matter of the validity of the available information sources, and probably the best the can be achieved is to approximate sufficiently close the true $\Phi(\pi)$ by $\tilde{\Phi}(\pi) := J(\pi, Q^*)$ given that the information sources led to the estimate $Q^*$ provided quite accurate predictions with respect to the real situation.

In order to formulate the problem under the framework of the model uncertainty, assume that the decision maker is provided with a collection of $n$ plausible probability models to play the role of the unknown model $Q$ by various agents of possibly varying validity. Let us denote by $\mathfrak{Q}:=\{Q_1,...,Q_n\}$ the set of the provided models and let also denote by $w=(w_1,...,w_n)\in\Delta^{n-1}$\footnote{If $x\in\Delta^{n-1}$ it means that $x_i\geq 0$ for all $i=1,2,...,n$ and $\sum_{i=1}^n x_i=1.$} the trustworthyness allocated by the decision maker to each one of the agents. According to the approach discussed in Section \ref{sec-2.2} the optimal probability model $Q_*$ should aggregate the information of the prior set $\mathfrak{Q}$ taking into account the importance of each opinion as determined by the decision maker, and of course the amount of aversion from the information provided that is allocated. If Kullback-Leibler divergence is employed as metric tool, then the set of admissible probability models can be written 
\begin{equation}
\mathcal{Q} := \{ Q\in\mathcal{P}\,\, :\,\, KL_w(Q, \mathcal{Q})\leq \eta \}
\end{equation}
where $\eta>0$ is set by the decision maker. Note here that if $\eta:=\min_Q KL_w(Q\|\mathcal{Q})$ then the decision maker is assumed to fully trust the provided prior models (no aversion) and the optimal probability model coincides with the KL-barycenter with respect to the weights $w$. On the other hand, if $\eta$ is chosen too large, then the decision maker does not appreciate at all the provided information and the optimal model $Q_*$ will be a very distorted version of the barycenter which may not provide realistic results. Clearly, it is very important the prior set to include at least one model that the decision maker could trust in order to avoid such behaviours.

Let us put some discussion on the model $Q$ and the reasoning of the model uncertainty framework in this case. The stochastic part of the problem under study is the random variable $Z = (Z_B, Z_R) \sim Q$ where $Q$ is the probability model under question. This model could be considered partially unknown or unknown as a whole. In the total ignorance case, everything needs to be retrieved by expert's opinion or espionage operations. Probability of successfull hits $p_B, p_R$ are considered unknown due to not sufficient knowledge of defensive and attacking potential and capabilities of the opponent's forces and as a result the dependence structure on exchanging fires cannot be also known. On the partial knowledge case, one may be sure about the attacking (or defensive) outcome of the opponent's side since the blue side's defensive (or attacking) capabilities are given. However, the dependence structure relies on various factors that is difficult to be precisely known, like battlefield and weather conditions, forces locations and movements, both sides weapon portfolios and tactics, etc. In fact, the dependence structure determines the manner that meet both sides their nominal attrition rates and the kind of correlation that characterizes them. If no dependence structure is assumed, the attrition rates are considered independent and the joint distribution binding could be described by the independence copula function. If the correlation structure is assumed as a linear one, then a Gaussian or Student-t copula maybe a good modeling approach. Of course there are much more choices, e.g. nonlinear dependencies can be represented by appropriate copula function choices. As a result, a plausible choice for the model $Q$ is crucial not only for the efficient modeling of states evolution $B$ and $R$, but also for the derivation of a realistic strategy $\pi$.

Following the discussion above the optimal strategy selection problem \eqref{opt-0} is reshaped to the robust control problem
\begin{eqnarray}\label{opt-1}
&&\max_{\pi \in [0,1]} \min_{Q \in \mathcal{Q}} J(\pi, Q) \\
&&\mbox{subject to}\nnb \\
&& \left\{ \begin{array}{ll}
B^{\pi}_{k \Delta t} = B^{\pi}_{(k-1) \Delta t} - r Z_{k,R} R^{\pi}_{(k-1)\Delta t} \Delta t, \,\,\, k=1,2,...,N, & B^{\pi}_0 = B_0 \\
R^{\pi}_{k \Delta t} = R^{\pi}_{(k-1) \Delta t} - \pi b Z_{k,B} B^{\pi}_{(k-1)\Delta t} \Delta t, \,\,\, k=1,2,...,N,& R^{\pi}_0=R_0\\
Z_k = (Z_{k,B}, Z_{k,R}) \sim Q,\,\,\, \forall k=1,2,...,N & 
\end{array} \right. \nnb
\end{eqnarray}
Treating first the inner problem, under the Kullback-Leibler divergence for any appropriate choice of the divergence parameter $\eta>0$, an explicit calculation of the optimal probability model $Q^*$ is possible according to \eqref{Qtheta} in terms of the probability density function. As a result the robust optimal strategy selection problem stated in \eqref{opt-1} is simplified to the initial formulation of problem \eqref{opt-0} where $Q^*$ is used in place of $Q$ condensing the available information of the prior set $\mathfrak{Q}$ and taking into account the decision maker's aversion preferences.

\subsubsection{A Numerical Recipe for the Optimal Strategy Selection}

After the aforementioned considerations and approaches, the final formulation of the under study optimal control (strategy) problem for given model risk preferences by the decision maker is the following:
\begin{eqnarray}\label{opt-2}
&&\max_{\pi \in [0,1]} J(\pi, Q^*) = \widetilde{\Phi}(\pi) \\
&&\mbox{subject to}\nnb \\
&& \left\{ \begin{array}{ll}
B^{\pi}_{k \Delta t} = B^{\pi}_{(k-1) \Delta t} - r Z_{k,R} R^{\pi}_{(k-1)\Delta t} \Delta t, \,\,\, k=1,2,...,N, & B^{\pi}_0 = B_0 \\
R^{\pi}_{k \Delta t} = R^{\pi}_{(k-1) \Delta t} - \pi b Z_{k,B} B^{\pi}_{(k-1)\Delta t} \Delta t, \,\,\, k=1,2,...,N,& R^{\pi}_0=R_0\\
Z_k = (Z_{k,B}, Z_{k,R}) \sim Q^*,\,\,\, \forall k=1,2,...,N & 
\end{array} \right. \nnb
\end{eqnarray}
where $\widetilde{\Phi}$ is defined as
\begin{equation}\label{objfun}
\widetilde{\Phi}(\pi) = \theta_1 \E_{Q^*}[ \varphi(\psi_1(X_T(\pi))) ] + \theta_2 \E_{Q^*}[ \varphi( \psi_2( X_T(\pi) ) ) ] + \theta_3 \psi( (1-\pi) B_0 )
\end{equation}
with 
\begin{eqnarray}
\psi_1( \pi; Z ) &=& \begin{bmatrix} 1 &-1 \end{bmatrix} \cdot X_{N \Delta t}(\pi; Z)\\
\psi_2( \pi; Z ) &=& \begin{bmatrix} 1 & \,\,\,\,\,\,0 \end{bmatrix} \cdot X_{N \Delta t}(\pi; Z) - B_{\min}
\end{eqnarray}
for any $Z = (Z_1, Z_2, ..., Z_N)'\sim Q^*$ and $X_{T}(\pi; Z)$ determined by equation \eqref{cf-sol} for $T := N \Delta t$ and can be written as
\begin{equation}
X_T(\pi; Z) = X_0 \cdot e^{\Delta t A_1(\pi; Z_1)} \cdot e^{\Delta t A_2(\pi; Z_2)} \cdots e^{\Delta t A_N(\pi; Z_N)}
\end{equation}
with the notation $A_k(\pi; Z_k)$ emphasizing on the direct dependence of the matrix on the parameter $\pi$ and its form on the realization of $Z_k$ for each $k=1,2,...,N$. The latter term causes some difficulties in taking first order optimality conditions with respect to $\pi$ since the differentiation of the matrix exponent is quite tricky and absolutely not straightforward. For this reason, for the derivation of both $\psi_1, \psi_2$ is followed the approach of complex step approximation proposed in \cite{Higham2010} which is computationally very simple and extremelly accurate. According to that, the directional derivative of $\psi_i$ is approximated as
\begin{equation}
\psi_i'(\pi) \approx Im \left( \frac{\psi_i(\pi + i h)}{h} \right) =: \hat{\psi}'_i, \,\,\, i=1,2,
\end{equation}
for very small $h>0$ and some direction $h$ where $Im(u)$ denotes the imaginary part of $u$.

A second problem in estimating the maximizer of \eqref{objfun} is that the current objective function depends on measure $Q^*$ and as a result an analytic expression for the objective function is computationally intractable. Based on that, a stochastic gradient scheme will be proposed, relying on sample trajectories from the combat dynamics in order to estimate both the value of the objective and its gradient with respect to $\pi$. In particular, the gradient of \eqref{objfun} is approximated through the function:
\begin{equation}
\widehat{\nabla}\widetilde{\Phi}(\pi) := \theta_1 \E_{Q^*}[ \varphi'(\psi_1(\pi; Z))\hat{\psi}'_1(\pi; Z) ] + \theta_2 \E_{Q^*}[ \varphi'(\psi_2(\pi; Z))\hat{\psi}'_2(\pi; Z) ] - \theta_3 B_0 \psi'((1-\pi)B_0)
\end{equation}

In this sense, a stochastic gradient scheme is described below for the approximation of the maximizer of the optimal control problem.

\begin{algorithm}[{\it Worst-Case Forces Optimal Allocation}]\label{alg-1}
Given the available initial forces $B_0, R_0$, a terminal time horizon $T>0$, a prior set $\mathfrak{Q}$ of models regarding the attrition dynamics, the decision maker's preferences $\theta, w, \eta>0$ repeat the following steps:
\begin{enumerate}
\item[{\bf 1.}] Determine the optimal probability model $Q_*$ as the minimizer of $\eqref{w-KLD}$ using the prior set $\mathfrak{Q}$ and the decision maker's preferences $(w,\eta)$.

\item[{\bf 2.}] Given the criteria weights $\theta$, determine the form of the objective function described in \eqref{objfun} with respect to the probability model $Q^*$. 

\item[{\bf 3.}] Approximate the maximizer $\pi_* = \arg\max_{\pi \in [0,1]} \widetilde{\Phi}(\pi)$ through the stochastic gradient ascent scheme
$$ \pi^{(\ell)} = \pi^{(\ell-1)} + \alpha_{\ell} \widehat{\nabla}\widetilde{\Phi}(\pi^{(\ell-1)}), \,\,\,\, \ell=1,2,...$$
where $\pi^{(0)} \in (0,1]$ an initial point, $\alpha_{\ell}>0$ the learning rate of the scheme at step $\ell$ and the iterations are terminated after certain convergence conditions are satisfied.  
\end{enumerate}
\end{algorithm}

\subsubsection{Extensions to the Mixed-Forces Case}

Although the discussion in the previous sections focused on the case where the strength of each side is measured as a whole without distinctions between different types of forces or weapons, the discussed framework can be very easily generalised to that case too. Assume that the blue side has a weapon portfolio of $m$ different systems while the red side strength is constituted by $q$ different systems. The blue side (decision maker) desires to select the best possible mixture of the available forces $\pi = (\pi_1, \pi_2, ..., \pi_m) \in [0,1]^m$ that will lead to the best result according to the commander's objectives. Under this perspective, the combat dynamical system may be represented (in a discretized version) as
\begin{eqnarray}\label{mixed-forces}
&& \left\{ \begin{array}{lll}
B^{\pi}_{j, k \Delta t} = B^{\pi}_{j, (k-1) \Delta t} - r_j Z^j_{k,R} G^{R}_j\left(R^{\pi}_{(k-1)\Delta t}\right) \Delta t, & B^{\pi}_{0,j} = B_{0,j}, & j=1,2,...,m, \\
R^{\pi}_{i, k \Delta t} = R^{\pi}_{i, (k-1) \Delta t} - \pi_i b_i Z^i_{k,B} G_i^{B}\left(B^{\pi}_{(k-1)\Delta t}\right) \Delta t, & R^{\pi}_{0,i}=R_{0,i}, & i=1,2,...,q, \\
Z_k = (Z_{k,B}, Z_{k,R}) \sim Q^*,& \forall k=1,2,...,N  &
\end{array} \right. 
\end{eqnarray}
where $G_i^{B}, G^{R}_j$ for $i=1,2,...,q$ and $j=1,2,...,m$ denote the attrition functions related to the attrition rate obtained by each weapon system of both sides. If the attrition functions are of linear form with respect to $B,R$, then a closed form representation can be obtained in the spirit of equation \eqref{cf-sol}. 

The choice of the dependence structure of the probability model $Q$ is very important here, since this choice should contain valuable information about the strategy that each side follows. For example, if a Gaussian copula is used to represent the dependence structure between the attrition rates, the correlation matrix parameters should quantify certain strategical segments of information about each side's positions, formations, correlations between the weapon systems that are used, etc. Although this information about the decision maker's side (the blue ones) can be considered as known, it does not hold the same for the opponents (red side) where different strategies may be used or diverging information about the side's stature are available. As a result, the proposed framework for treating model uncertainty in the simplified version of this problem, seems to be a very appropriate and effective tool in choosing a plausible probability model $Q^*$ characterizing the stochastic dynamics of the combat (after appropriate weighting and validation of the incoming imformation by the various experts/agents). 

As a result, this more general version of weapon resources allocation problem can be stated as a minor variation of the problem \eqref{opt-2} while the numerical approach described in algorithm \ref{alg-1} can be also easily modified for approximating the problem's maximizer.

\section{Conclusion}

In this work was made a first attempt in aggregating the optimal decision making theory with the model uncertainty framework on the field of optimal strategy selection in combat situations. A stochastic variation of the famous Lanchester aimed-fire model was used as a basis for modeling the combat dynamics and a robust decision rule was constructed and proposed under the framework of model ambiguity describing the dynamics of the state variables (combat strength) and taking into account the decision maker's preferences and expected goals. Upon appropriate discretizations, a stochastic gradient method has been proposed for approximating the maximizer of the robust control strategy selection problem.

\section*{Acknowledgements}

The author would like to express his gratitude to Professor S. Kyritsi-Yiallourou (HNA) for introducing him to the exciting world of warfare modeling and for the very interesting and fruitful
discussions on the matter. He is also extremely thankful to his good friend and colleague Professor
A. N. Yannacopoulos (AUEB) for initiating him in the art of robust decision making at theoretical
and empirical level.


\end{document}